\documentclass[12pt,intlimits,namelimits]{amsart}
\usepackage{amssymb}
\begin{document}

\vspace*{10mm}

\vspace{5mm}

\title[Ideals in some Rings of Nevanlinna--Smirnov 
Type]{Ideals in some Rings of Nevanlinna--Smirnov Type}
\author[]{Romeo Me\v{s}trovi\'c}
\address{University of Montenegro, Maritime Faculty, \endgraf
85330 Kotor, Montenegro}
   \begin{abstract}
Let $N^p$ $(1<p<\infty)$ denote the algebra of holomorphic functions
in the open unit disk, introduced by I.~I.~Privalov with the notation
$A_q$ in [8]. Since $N^p$ becomes a ring of Nevanlinna--Smirnov type
in the sense of Mortini [7], the results from [7] can be applied to
the ideal structure of the ring $N^p$. In particular, we observe that
$N^p$ has the Corona Property. Finally, we prove the $N^p$-analogue of the
Theorem 6 in [7], which gives sufficient conditions for an
ideal in $N^p$, generated by a finite number of inner functions, to be equal
to the whole algebra $N^p$.
\end{abstract}
\maketitle

\section{Introduction and Preliminaries}

Let $D$ denote the open unit disk in the complex plane and let $T$ denote
the boundary of $D$. Let $L^p(T)$ $(0<p\leqslant\infty)$ 
be the familiar Lebesgue spaces on $T$. The \emph{Nevanlinna class} $N$ 
is the set of all functions $f$
holomorphic on $D$ such that
\[
\sup_{0<r<1}\int_0^{2\pi}\log^+
\big\lvert
f\big(re^{i\theta}\big)\big\rvert\frac{\mathrm{d}\theta}{2\pi}
<\infty{,}
\]
where $\log^+|x|=\max(\log|x|,0)$.

The \emph{Smirnov class} $N^+$ consists of those functions $f\in N$ for which
\[
\lim_{r\to 1}\int_0^{2\pi}\log^+
\big\lvert
f\big(re^{i\theta}\big)\big\rvert\frac{\mathrm{d}\theta}{2\pi}=
\int_0^{2\pi}\log^+
\big\lvert
f^*\big(e^{i\theta}\big)\big\rvert\frac{\mathrm{d}\theta}{2\pi}
<\infty{,}
\]
where $f^*$ is the boundary function of $f$ on $T$, i.e.,
$$
f^*\big(e^{i\theta}\big)=\lim_{r\to 1}f\big(re^{i\theta}\big)
$$
is the radial limit of $f$ which exists for almost every $e^{i\theta}$.

Recall that the \emph{Hardy space} $H^p$ $(0<p\leqslant \infty)$ consists of
all functions $f$ holomorphic in $D$, which satisfy
\[
\sup_{0< r<1}\int_0^{2\pi}\big\lvert
f\big(re^{i\theta}\big)\big\rvert^p\frac{\mathrm{d}\theta}{2\pi}
<\infty.
\]
If $0<p<\infty$, and which are bounded when $p=\infty$;
\[
\sup_{z\in D}|f(z)|<\infty.
\]

Following R. Mortini [7], a ring $R$ satisfying $H^\infty\subset R\subset N$
is said to be of \emph{Nevanlinna--Smirnov type} if every function $f\in  R$
can be written in the form $g/h$, where $g$ and $h$ belong to $H^\infty$ and $h$ is
invertible element in $R$. This is true of $N$ itself and the Smirnov class
$N^+$; hence the name (see [1, Chapter 2]). Further,
Mortini noted that by a result of M. Stoll [10], the space $F^+$, the containing
Fr\'echet envelope for $N^+$, consists of those functions $f$ holomorphic in $D$
satisfying
\[
\limsup_{r\to 1}(1-r)\log M(r,f)=0
\]
with $M(r,f)=\max_{|z|=r}|f(z)|$ (see Yanagihara [11]).

Namely, Stoll [10] proved that
$F^+\cap N=\{f/S_\mu:f\in N^+$, $S_\mu$ is a singular inner function with $\mu$
a nonnegative continuous singular measure\}.

The class $N^p$ $(1<p<\infty)$ consists of all holomorphic functions $f$ on $D$
for which
\[
\sup_{0<r<1}\int_0^{2\pi}\Big(\log^+
\big\lvert
f\big(re^{i\theta}\big)\big\rvert\Big)^p\frac{\mathrm{d}\theta}{2\pi}
<\infty{.}
\]
These classes were introduced in the first edition of Privalov's book
[8, p.~93], where $N^p$ is denoted as $A_q$. It is known [6] that
\[
N^q\subset N^p\; (q>p),\quad \mathop{\bigcup}_{p>0} H^p\subset \mathop{\bigcap}_{p>1}N^p,
\quad\text{and}\quad   \mathop{\bigcup}_{p>1} N^p\subset N^+,
\]
where the above containment relations are proper.

\vspace{3mm}

\noindent{\bf Theorem A}
([8, p.~98]).
{\it A function $f\in N^p\setminus\{0\}$ has a unique factorization of the form
\[
f(z)=B(z)S(z)F(z)
\]
where $B(z)$ is the Blaschke product with respect to zeros
$\{z_k\}\subset D$ of $f(z)$, $S(z)$ is a singular inner function, $F(z)$
is an outer function in $N^p$, i.e.,
\[
B(z)=z^m\prod_{k=1}^\infty\frac{|z_k|}{z_k}\,
\frac{z_k-z}{1-\bar{z}_kz}
\]
with $\sum_{k=1}^\infty\big(1-|z_k|\big)<\infty$, $m$ a nonnegative integer,
\[
S(z)=\exp\left(-\int_0^{2\pi}\frac{e^{it}+z}{e^{it}-z}\,
\mathrm{d}\mu(t)\right)
\]
with positive singular measure $\mathrm{d}\mu$, and
\[
F(z)=\lambda\exp\left(\frac{1}{2\pi}\int_0^{2\pi}\frac{e^{it}+z}{e^{it}-z}
\log\big\lvert f^*\big(e^{it}\big)\big\rvert\,\mathrm{d}t\right),
\]
where $|\lambda|=1$,  $\log|f^*|$ and $\big(\log^+|f^*|\big)^p$ belong
to $L^1(T)$.

Conversely, every such product $B(z)S(z)F(z)$ belongs to $N^p$.}

\vspace{3mm}

\noindent{\it Remark.} If we exclude only the condition $\big(\log^+|f^*|\big)^p\in L^1(T)$
from Theorem A, we obtain the well known canonical factorization theorem
for the class $N^+$ (see [1, p.~26] or [8, p.~89]). Recall that a function of the
form $\varphi(z)=B(z)S(z)$ is called an \emph{inner function}, while the
above function $F(z)$ for which $\log|f^*|\in L^1(T)$ is called an
\emph{outer function}.

\vspace{3mm}

By Theorem A, it is easy to show (see [2], where
$N^p$ is denoted as $N^+_\alpha$) that a function $f$ is in $N^p$ if and
only if it can be expressed as the ratio $g/h$, where $g$ and $h$ are in
$H^\infty$, and $h$ is an outer function with
$\log|h^*|\in L^p(T)$. Clearly, such function $h$ is an invertible
element of $N^p$, and hence we have the following result.

\vspace{3mm}

\noindent{\bf Theorem B.}
{\it $N^p$ $(1<p<\infty)$ is a ring of Nevanlinna--Smirnov type.}

\vspace{3mm}

In Section 2, the ideal structure of subrings $N^p$ of $N$ is described as
consequences of the results in [7, Sections 1 and 3] given for an arbitrary
ring of Nevanlinna--Smirnov type. In the next section, we note that the algebra
$N^p$ has the Corona Property. Finally, we prove two theorems which generalize
the results from [7, Satz 5 and Satz 6] obtained for the classes $H^1$ and $N^+$.

\section{Ideals in $N^p$}

In this Section, as an application of Theorems A and B and the results of
Mortini in [7], we obtain some facts about the ideal structure of the
algebra $N^p$. We say that an ideal $I$ in $H^\infty$ is the trace of an
ideal $J$ in $N^p$ if $I=J\cap H^\infty$. The following result is an
immediate consequence of Theorems A, B and [7, Satz 1, Satz 2].

\vspace{3mm}

\noindent{\bf Theorem 1.}
{\it
 An ideal $I$ in $H^\infty$ is the trace of an ideal $J$ in $N^p$ if
 and only if the following condition is satisfied: If $f\in I$, $F$
 is an outer function with $\log|F^*|\in L^p(T)$, and if $fF\in
 H^\infty$, then $fF\in I$.
 In this case, $J$ is a unique ideal in $N^p$ with $I=J\cap H^\infty$,
 and there holds $J=IN^p$.}

\vspace{3mm}

Further, by the above theorem, it follows immediately the following theorem.

\vspace{3mm}

\noindent{\bf Theorem 2.}
{\it
Suppose that $I$ is an ideal in $H^\infty$ such that $f\in I$ implies
that the inner factor of $f$ also belongs to $I$.
Then $I$ is the trace of an ideal $J$,  in $N^p$ and there holds $J=IN^p$.}

\vspace{3mm}

\noindent{\it Remark.}
It remains an open question is it true the converse of Theorem 2.
While this is true for the Nevanlinna class and the Smirnov class
[7, Korrolar 1 and Korrolar 2, resp.], here the corresponding problem
is complicated by the fact that there exist outer functions which are not
invertible in $N^p$. For example, the converse of Theorem 2 holds if we
suppose in addition that $I$ is a closed subset of $N^p$. Namely, by
Theorem 2 of [6] and the condition from Theorem 1, it follows that
whenever $f$ is in $I$, then necessarily the inner factor of $f$ is also in $I$.

\vspace{3mm}

Recall that an ideal $P$ in a ring $R$ is prime if whenever $fg\in P$,
$f,g\in R$,
then either $f$ or $g$ is in $P$. Using the characterization of the invertible
elements in $N^p$, by [7, Satz 3], we obtain the following.

\vspace{3mm}
\noindent{\bf Theorem 3.}
{\it
A prime ideal $P$ in $H^\infty$ is the trace of some prime ideal $Q$ in
$N^p$ if and only if $P$ contains no outer functions $F$ for which
$\log|F^*|\in L^p(T)$. When this is the case, $Q$ is a unique prime ideal in
$N^p$ with this property, and there holds $Q=PN^p$.}

\vspace{3mm}

\noindent{\it Remark.}
By [6, Theorem 3], every prime ideal of $N^p$ which is not dense in $N^p$ is
equal to the set of functions in $N^p$ vanishing at a specific point of $D$.
The analogous result for the class $N^+$ is proved in [9, Theorem 1].

\vspace{3mm}

An ideal $J$ in the ring $R$, $H^\infty\subset R\subset N$, is called finitely
generated if there exist elements $f_1,\ldots,f_n\in R$ such that
\[
J=(f_1,\ldots,f_n)=\left\{\sum_{i=1}^n g_if_i:g_i\in R\right\}{.}
\]
\vfill\eject

If $n$ can be chosen to be one, then $J$ is a principal ideal.
A ring $R$ is said to be coherent if the intersection of two finitely
generated ideals is finitely generated. Using the result in [5] that $H^\infty$
is a coherent ring, it is shown in [7, Satz 7] that this is true of all
rings of Nevanlinna--Smirnov type. In particular, we have the following

\vspace{3mm}

\noindent{\bf Theorem 4.}
{\it
$N^p$ is a coherent ring for all $p>1$.}

\section{The Corona Property}

We say that a commutative ring $R$ with unit of holomorphic functions
on $D$ has the \emph{Corona Property} if the ideal generated by
$f_1,\ldots,\allowbreak f_n\in R$ is equal to $R$ if and only if there is an
invertible element $f$ of $R$ such that
\[
|f(z)|\leqslant\sum_{i=1}^n |f_i(z)|\quad (z\in D).
\]
This definition is motivated by the famous Corona Theorem of Carleson
(for example see [3, p.~324], or [1, p.~202]), which states that the algebra
$H^\infty$ of all bounded holomorphic functions on $D$ has the Corona Property.
Mortini noted [7, Satz 4] that by a result of Wolff [3, p.~329], it is easy to
show that every ring of Nevanlinna--Smirnov type has the Corona Property. Hence,
we have the following theorem.

\vspace{3mm}

\noindent{\bf Theorem 5.}
{\it
The algebra $N^p$ has the Corona Property.}

\vspace{3mm}

\noindent{\it Remark.}
It is proved in [4, Theorem 7] that there exists a subalgebra of $N$
containing $N^+$ without the Corona Property.

\vspace{3mm}

A sequence $\{z_k\}\subset D$ is called an \emph{interpolating sequence}
(for $H^\infty$) if for every bounded sequence $\{\omega_k\}$ of
complex numbers there exists a function $f$ in $H^\infty$ such that
$f(z_k)=\omega_k$ for every $k$. An \emph{interpolating Blaschke
product} is a Blaschke product whose (simple) zeros form an
interpolating sequence.

The following two theorems generalize Theorems 5 and 6 in [7], respectively. We follow [7]
for the proofs of theorems bellow.

\vspace{3mm}

\noindent{\bf Theorem 6.}
{\it
Let $0<q<\infty$ and let $I=(f_1,\ldots,f_n)$ be a finitely generated
ideal in $H^\infty$. Assume that the ideal $I$ contains a zero-free
holomorphic function $F$ on $D$ such that $\log F \in H^q$. Then}
\[
\sum_{k=1}^\infty\big(1-|z_k|^2\big)\left|\log \left(|f_1(z_k)|+\cdots
+|f_n(z_k)|\right)\right|^q<\infty{.}
\]

\begin{proof}
Let $g_1,\ldots,g_n$ be functions in $H^\infty$ such that
$F=\sum_{k=1}^n f_kg_k$. Then there exists a positive constant $C_1$
such that
\[
|F(z)|\leqslant C_1 \sum_{k=1}^n |f_k(z)|\quad \text{for all}\quad z\in D.
\]
Put $S(z)=\sum_{k=1}^n |f_k(z)|$, and suppose that $S(z)\leqslant C_2$ for
all $z\in D$, with a positive constant $C_2$. Using the fact that
$\sum_{k=1}^\infty (1-|z_k|)=C_3<\infty$, applying the inequality
$(a+b)^q\leqslant C_4\big(a^q+b^q\big)$ with $C_4=2^{\max(q,1)-1}$,
$a,b\geqslant 0$, and the main interpolation theorem for the class
$H^q$ [1, p.~149], we obtain
\[
\begin{split}
&\sum_{k=1}^\infty\big(1-|z_k|^2\big)\big|\log S(z_k)\big|^q\\
&=\sum_{k=1}^\infty\big(1-|z_k|^2\big)\left(\log^+ S(z_k)+\log^+\frac{1}{S(z_k)}\right)^q\\
&\leqslant \sum_{k=1}^\infty\big(1-|z_k|^2\big) C_4
\left(\left(\log^+ S(z_k)\right)^q+\left(\log^+\frac{1}{S(z_k)}\right)^q\right)\\
&\leqslant C_4(C_2)^q\sum_{k=1}^\infty\big(1-|z_k|^2\big)+C_4 \sum_{k=1}^\infty\big(1-|z_k|^2\big)
\!\!\left(\log^+\frac{C_1}{|F(z_k)|}\right)^q\\
&\leqslant 2C_3C_4(C_2)^q+C_4 \sum_{k=1}^\infty\big(1-|z_k|^2\big)
\left|\log\frac{F(z_k)}{C_1}\right|^q<\infty{.}
\end{split}
\]
This gives the desired result.
\end{proof}

\noindent{\bf Theorem 7.}
{\it
Assume that $I$ is an ideal in $N^p$ generated by inner functions
$\varphi_1,\ldots,\varphi_n$, and suppose that $I$ contains an interpolating
Blaschke product $B$ with zeros $\{z_k\}$ such that
\[
\sum_{k=1}^\infty\big(1-|z_k|^2\big)\left|\log\left(|\varphi_1(z_k)|+\cdots
|\varphi_n(z_k)|\right)\right|^p<\infty{.}
\]
Then $I=N^p$.}

\begin{proof}
Put $c_k=\sum_{i=1}^n |\varphi_i(z_k)|^2$ for all $k$. Using the inequalities
\[
\frac{\left( \sum\limits_{i=1}^n |\varphi_i(z_k)|\right)^2}{n}\leqslant c_k\leqslant
\left( \sum_{i=1}^n |\varphi_i(z_k)|\right)^2,
\]
it is routine to estimate that
\[
|\log c_k|\leqslant 2\left|\log\left( \sum_{i=1}^n|\varphi_i(z_k)|\right) \right|
+\log n,
\]
whence by the inequality $(a+b)^p\leqslant 2^{p-1}(a^p+b^p)$,
$a,b\geqslant 0$,
using the assumption of the theorem, we have
\[
\begin{split}
\qquad &\sum_{k=1}^\infty\big(1-|z_k|^2\big)|\log c_k|^p\\
&\leqslant 2^{2p-1}\sum_{k=1}^\infty \big(1-|z_k|^2\big)\left|\log
\left(\sum_{i=1}^n\big|\varphi_i(z_k)\big|\right)\right|^p \\
&+2^{p-1}\log^p n \sum_{k=1}^\infty \big(1-|z_k|^2\big)<\infty{.}
\end{split}
\]
Hence, by the theorem of Shapiro and Shields [1, p.~149, Theorem
9.1], there exists a function $g\in H^p$ with $g(z_k)=\log c_k$ for
every $k$. It is easy to verify that the function $F=\exp g$ is
invertible in $N^p$, and there holds $F(z_k)=c_k$.

The rest of the proof is the same as that of [7, Satz 6]. To complete
the proof, we write this part.

Since $\{z_k\}$ is an interpolating sequence, by [1, p.~149], we
know that there exist functions $f_i\in H^\infty$ $(i=1,\ldots,n)$
such that for every $i=1,\ldots,n$ there holds
\[
f_i(z_k)=\overline{\varphi_i(z_k)}\quad (k=1,2,\ldots)
\]
Observe that the function $F-\sum_{i=1}^n f_i\varphi_i$ is in $N^p$,
and that there holds
\[
F(z_k)-\sum_{i=1}^n f_i(z_k)\varphi_i(z_k) =c_k-
\sum_{i=1}^n |\varphi_i(z_k)|^2=0\quad (k=1,2,\ldots)
\]
Hence, by Theorem A, there exists a function $h\in N^p$ such that
\[
F-\sum_{i=1}^n f_i\varphi_i=Bh{.}
\]
This shows that $F$ belongs to the ideal $(\varphi_1,\ldots,\varphi_n,B)=I$.
Since $F$ is an invertible element in $N^p$, it follows that
\[
I=(\varphi_1,\ldots,\varphi_n)=N^p{.}
\]
This completes the proof of the theorem.
\end{proof}

\end{document}